\documentclass{amsart}
\usepackage{amsmath,amsfonts,amsthm}

\newtheorem{Proposition}[equation]{Proposition}
\newtheorem{Lemma}[equation]{Lemma}
\newtheorem{Theorem}[equation]{Theorem}

\newtheorem{Definition}[equation]{Definition}
\newtheorem{Remark}[equation]{Remark}

\newtheorem{Example}[equation]{Example}

\def\S{\mathcal{S}}
\def\C{\mathcal{C}}
\def\gS{\operatorname{gr} \mathcal{S}^f}
\def\Z{\mathbb{Z}_2}
\def\eZ#1{(\mathbb{Z}_2^{#1})^{\text{ev}}}
\def\oZ#1{(\mathbb{Z}_2^{#1})^{\text{odd}}}
\def\M{\mathcal{M}}
\def\P{\mathcal{P}}
\def\eM{\mathcal{M}^\text{ev}}
\def\eP{\mathcal{P}^\text{ev}}
\def\deg{\operatorname{deg}}
\def\x{\hat{x} }
\def\c{\hat{c} }
\def\s{\hat{s} }

\def\Acy{{\stackrel{\rightarrow}{A}}}
\def\acy{{\stackrel{\rightarrow}{\alpha}}}
\def\sp{\operatorname{supp}}

\def\Cen{\operatorname{Cent}}

\def\gr{\operatorname{gr}}

\begin{document}

\title{Centers of cyclotomic Sergeev superalgebras}
\author{Oliver Ruff}

\address{Department of Mathematics, University of Toledo, Toledo, OH 43606, USA}
\email{oliver.ruff@utoledo.edu}

\begin{abstract}
We prove that the natural map from the center of the affine Sergeev superalgebra to the even center of any cyclotomic Sergeev superalgebra of odd level is surjective, hence that the even center of a cyclotomic Sergeev superalgebra of odd level consists of symmetric functions in the squares of its polynomial generators.  
\end{abstract}

\maketitle

\section{Introduction}

This paper is an attempt at a twisted version of Brundan's description of the centers of the type A degenerate cyclotomic Hecke algebras $H^f_d$, and thereby his classification of their blocks. These algebras arise as finite-dimensional quotients of the corresponding degenerate affine Hecke algebra $H_d$, whose center is well-known to consist of all symmetric polynomials in its polynomial generators, so it is easy to write down central elements of $H^f_d$ corresponding to the elements of $Z(H_d)$. However, it is far from obvious that these elements constitute $Z(H^f_d)$ in its entirety.

In \cite{B}, Brundan filters $H^f_d$ by polynomial degree and studies the associated graded object, which is a twisted tensor product of a truncated polynomial ring with the group algebra of the symmetric group. By examining the centralizer of the polynomial subalgebra, and finding within it the fixed points of a natural conjugation action by the symmetric group, an explicit description of the center of the graded algebra is obtained in terms of certain elements called \emph{colored cycles}. This provides the needed upper bound on the dimension of $Z(H^f_d)$.

In place of $H_d$, we study the \emph{affine Sergeev superalgebra} $\S_d$ introduced by Nazarov in \cite{N}. The appropriate finite-dimensional quotients are the \emph{cyclotomic Sergeev superalgebras} $\S^f_d$, where $f$ is a monic polynomial of degree $l$ having a certain specified form \cite[section 3-e]{BK}. (Full definitions follow in Section~\ref{prelim}.) The center of $\S_d$ is also well-known to be generated by the squares of the polynomial generators, and so the logical question is whether $Z(\S^f_d)$ is equal to the image of $Z(\S_d)$ under the quotient map. We prove the following:

\begin{Theorem}
Let $l$ be odd. Then $Z(\S^f_d)$ is the free $R$-module consisting of all symmetric polynomials in the squares $\x_1^2,\dots,\x_d^2$ of its polynomial generators.
\end{Theorem}

When $l$ is even we show that the center of the graded superalgebra is strictly larger than the image of $Z(\S_d)$ in $Z(\S^f_d)$. However, by itself this is not sufficient to conclude anything about the rank of $Z(\S^f_d)$ in that case.

This result provides the expected classification of the blocks of $\S^f_d$ for odd $l$, which partially corrects a gap in the existing literature. Specifically, this is the degenerate version of the claim made in \cite[Corollary 8.13]{BK}. 

\emph{Acknowledgements:} I am grateful to Jonathan Brundan and Alexander Kleshchev for helpful comments and advice. This work was begun while visiting the University of Sydney in 2007, and I am indebted to Andrew Mathas and the algebra group there for their hospitality.

\section{Preliminaries}\label{prelim}

\subsection{Notation}\label{notation} 

Fix $l,d \geq 1$ and a commutative ring $R$. For any $a\geq 1$, write $I_a$ for the index set $\{1,\dots,a\}$.

For $A= \{i_1,\dots,i_a \} \subseteq I_d$ and $\alpha \in \Z^a$, write $\sp(\alpha,A)$ for the set $\{ \ i_j \ | \ \alpha_j = 1 \ \}$. Define $|\alpha| := \sum_{i=1}^a \alpha_i = |\sp(\alpha,I_a)| $, and call $\alpha$ even or odd according to the parity of $|\alpha|$. Write $\eZ{a}$ and $\oZ{a}$ for the set of all even and odd elements of $\Z^a$, respectively. 

For $i \in I_a$, write $1_i$ for the element of $\Z^a$ with a $1$ in the $i$th position and $0$s elsewhere. 

For $1 < j \in I_a$ and $\alpha \in \Z^a$, define 
$$\alpha^{(j)} := \alpha + 1_j + 1_{j-1}$$
and
$$\alpha^{(1)} := \alpha + 1_1 + 1_a .$$

\subsection{(Super)algebras}

Let $R_l[x_1,\dots,x_d]$ denote the truncated polynomial algebra $$ R[x_1,\dots,x_d]/(x_1^l,\dots,x_d^l) .$$ 
Let $R \Sigma_d$ be the group algebra of the symmetric group. We regard both of these as superalgebras concentrated in degree $0$. (The word ``super'' will often be suppressed in what follows.)

Write $C_d$ for the Clifford algebra on odd generators $c_1,\dots,c_d$, with relations
$$ c_i^2 = 1,  \hspace{.5in} c_i c_j = - c_j c_i \text{ for } i \neq j \in I_d .$$

For $A = \{i_1,\dots,i_a \} \subseteq I_d$, write $\sigma_A$ for the cycle $(i_1 \ \dots \ i_a) \in R \Sigma_d$. Write $R_l[A]$ for the algebra $R_l[x_i \ | \ i \in A]$ and $C[A]$ for the Clifford algebra generated by $\{c_i \ | \ i \in A \}$. 

\begin{Definition}{\rm 
Let $\S_d$ be the \emph{affine Sergeev superalgebra}. As an $R$-module, it is free on Coxeter generators $\s_1,\dots,\s_{d-1}$, polynomial generatora $\x_1,\dots,\x_d$, and Clifford generators $\c_1,\dots,\c_d$, which are then subject to the following relations:
\begin{enumerate}
\item[$\bullet$] Coxeter generators satisfy the relations in $R\Sigma_d$; polynomial generators commute with one another; Clifford generators satisfy the relations in $C_d$.
\item[$\bullet$] $\x_i \c_i = - \c_i \x_i$, $\x_i \c_j = \c_j \x_i$ for $i, j \in I_d$ with $j \neq i$.
\item[$\bullet$] $\s_i \c_i = \c_{i+1} \s_i$, $\s_i \c_j = \c_j \s_i$ for $i,j \in I_{d-1}$ with $j \neq i, i+1$.
\item[$\bullet$] $\s_i \x_i = \x_{i+1} \s_i - 1 - \c_i \c_{i+1}$, $\s_i \x_j = \x_j \s_i$ for $i \in I_{d-1}$, $j \in I_d$ with $j \neq i, i+1$.
\end{enumerate} 
}\end{Definition}

We record the following identity in $\S_d$, for $i \in I_{d-1}$ and $n \geq 1$:
\begin{equation*}\label{affpowerrel} \s_i \x_i^n = \x_{i+1}^n \s_i - \sum_{j=0}^{n-1} \x_i^j \x_{i+1}^{n-1-j} (1 + (-1)^j \c_i \c_{i+1}) . \end{equation*}

The following result is well-known:

\begin{Theorem}\cite{K}\label{affcenter}
$Z(\S_d)$ is the set consisting of all symmetric polynomials in the $\x_1^2, \dots, \x_d^2$.  
\end{Theorem}

\begin{Definition}\label{fdef}{\rm
Let $f(x) = x^l + b_{l-2}x^{l-2} + b_{l-4} x^{l-4} + \dots \in R[x]$ be a monic polynomial with the property that the terms appearing in $f$ all have either even or odd degree. Then the \emph{cyclotomic Sergeev superalgebra} $\S^f_d$ is the quotient of $\S_d$ by the two-sided ideal generated by $f(\x_1)$. The parameter $l$ is called the \emph{level} of $\S^f_d$.
}\end{Definition}

As a consequence of Theorem~\ref{affcenter}, the image of every symmetric polynomial in the $\x_1^2,\dots,\x_d^2$ is central in $S^f_d$. The purpose of the rest of the paper is to prove that that there are no other central elements. 

We give a filtration $ F_0 \S_d^f \subset F_1 \S_d^f \subset \dots $ of $\S_d^f$ by declaring Coxeter and Clifford generators to lie in filtered degree $0$, and polynomial generators to lie in filtered degree $1$. Write $\gS_d$ for the associated graded object. By the PBW theorem for cyclotomic Sergeev superalgebras \cite[Theorem 15.4.6]{K}) there is an isomorphism of superalgebras
$$ \gS_d \cong (R_l[x_1,\dots,x_d] \hat{\otimes} R \Sigma_d) \otimes C_d $$
where $\hat{\otimes}$ is the usual twisted tensor product. We identify $R_l[x_1,\dots,x_d]$, $R \Sigma_d$, and $C_d$ with the corresponding subalgebras of $\gS_d$ when convenient. 

\begin{Remark}{\rm
We will abuse notation by writing $\x_i$, $\c_i$, and $\s_i$ for the generators of both $\S_d$ and $\S^f_d$. We write $x_i$, $c_i$ and $s_i$ for the generators of $\gS_d$, where we will be working from now until section~\ref{penultimate}.
}\end{Remark}

For $ A \subseteq I_d$ and $z \in \gS_d$, we say $z$ is of \emph{maximal degree} with respect to $A$ if $x_i^{l-1}$ divides $z$ for all $i \in A$. An element of the form $f \sigma_A c \in \gS_d$ with $f \in R_l[A]$ and $c \in C[A]$ is called a \emph{maximal degree cycle} if it is of maximal degree with respect to $A$. Define the \emph{length} of an element $f \sigma c \in \gS_d$ with $f \in R_l[d]$, $\sigma \in \Sigma_d$, and $c \in C_d$ to be the length of $\sigma$.

\subsection{Signs}

\begin{Definition}{\rm For $\alpha \in \Z^a$, define an associated $a$-tuple of signs 
$$ \epsilon^\alpha = (\epsilon^\alpha_1,\dots,\epsilon^\alpha_a) \in \{ \pm 1 \}^a$$ 
by
$$ \epsilon^\alpha_i := \prod_{j < i } (-1)^{\alpha_j} $$
for $i \in I_a$, where we understand $\epsilon^\alpha_1 = 1$ for all $\alpha$. 
}\end{Definition}

\begin{Remark}\label{signprop}{\rm
\begin{enumerate}
\item[(i)] $\alpha \in \Z^a$ lies in $\eZ{a}$ if and only if
$$ \epsilon^\alpha_a = (-1)^{\alpha_a} .$$
\item[(ii)] For $A = \{i_1,\dots,i_a \} \subseteq I_d$, $\alpha \in \Z^a$, and $j \in I_a$ we have 
\begin{eqnarray*} c_{i_j} c_\alpha (A) & = & \epsilon^\alpha_j c_{\alpha + 1_j} (A) , \hspace{.5in} \text{and} \\
c_\alpha (A) c_{i_j} & = & (-1)^{\alpha_{j+1}+\dots+\alpha_d} c_{\alpha + 1_j} (A) . \end{eqnarray*}
Moreover, if $\alpha \in \eZ{a}$, then this last expression equals $\epsilon^\alpha_{j+1} c_{\alpha+1_j}(A)$. (Here the sign is to be interpreted as $1$ if $j=a$.)
\item[(iii)] For $1 < j \in I_a$, 
$$ \epsilon^{\alpha^{(j)}}_i = \left\{ \begin{array}{ccl} 
       \epsilon^\alpha_i  & & i \neq j \\
       &&\\
       - \epsilon^\alpha_i & & i = j, \end{array} \right. $$
whereas
$$ \epsilon^{\alpha^{(1)}}_i = \left\{ \begin{array}{ccl} 
      - \epsilon^{\alpha}_i  & & i \neq 1 \\
       &&\\
       \epsilon^\alpha_1 & & i = 1. \end{array} \right. $$  
\end{enumerate}
}\end{Remark}

\subsection{Polynomials}

Most of our results ultimately depend upon the good behaviour of a certain class of polynomials, which are generalizations of the ones defined in chapter 2 of \cite{B}. Here we define them and establish some of their properties.

\begin{Definition}{\rm
For $A = \{ i_1,\dots,i_a \} \subseteq I_d$, $r \geq 0$, and $\alpha \in \eZ{a}$, the polynomial $h_r^\alpha(A)$ is defined as follows:
$$ h_r^\alpha(A) = \sum_{r_1 + \dots + r_a = (a-1)(l-1)+r } (\epsilon^\alpha_1 x_{i_1})^{r_1} \dots (\epsilon^\alpha_{i_a} x_{i_a})^{r_a} .$$
}\end{Definition}

\begin{Lemma}\label{commhc}{\rm
\begin{enumerate}
\item[(i)] For $1 < j \leq d$, $c_{i_j} \ h_r^\alpha(A) \ = \ h_r^{\alpha^{(j)}}(A) \ c_{i_j}$.
\item[(ii)] $c_{i_1} \ h_r^\alpha(A) \ = \ (-1)^{(a-1)(l-1)+r}  h_r^{\alpha^{(1)}}(A) \ c_{i_1}$.
\end{enumerate}
}\end{Lemma}

\begin{proof}
\begin{enumerate}
\item[(i)] Follows immediately from Remark~\ref{signprop}~(iii).
\item[(ii)] \begin{eqnarray*}
c_{i_1} \  h^\alpha_r (A) & = & c_{i_1} \sum_{r_1 + \dots + r_a = (a-1)(l-1)+r } (\epsilon^\alpha_1 x_{i_1})^{r_1} \dots (\epsilon^\alpha_{i_a} x_{i_a})^{r_a} \\
& = & \sum_{r_1,\dots,r_a} \Big( (-\epsilon^\alpha_1 x_{i_1} )^{r_1} \prod_{1<j \in I_a} (\epsilon^\alpha_j x_{i_j} )^{r_j} \Big) c_{i_1} \\
& = & \sum_{r_1,\dots,r_a} \Big( (-1)^{r_1 +\dots+r_a} (\epsilon^\alpha_1 x_{i_1})^{r_1} \prod_{1<j \in I_a} ( - \epsilon^\alpha_j x_{i_j})^{r_j} \Big) c_{i_1} \\
& = & (-1)^{(a-1)(l-1)+r} h^{\alpha^{(1)}}_r (A) \ c_{i_1} \hspace{.5in} \text{by Remark~\ref{signprop}~(iii).} \end{eqnarray*} \end{enumerate}  \end{proof}

\begin{Lemma}\label{commhx} For $A =\{i_1,\dots,i_a \}$, $\alpha \in \eZ{a}$, $j \in I_a$, and $r \geq 0$, we have 
\begin{equation}\label{commhxform} x_{i_j} h^\alpha_r (A) = (-1)^{\alpha_j} h^\alpha_r(A) x_{\sigma_A(i_j)}. \end{equation} \end{Lemma}

\begin{proof}
We have
\begin{eqnarray*}
x_{i_j} h^\alpha_r(A) & = & x_{i_j} \sum_{r_1+\dots+r_a=(a-1)(l-1)+r} (\epsilon^\alpha_1 x_{i_1})^{r_1} \dots (\epsilon^\alpha_a x_{i_a})^{r_a} \\
& = & \epsilon^\alpha_j \ \sum_{r_1,\dots,r_a} (\epsilon^\alpha_1 x_{i_1})^{r_1} \dots (\epsilon^\alpha_j x_{i_j})^{r_j+1} \dots (\epsilon^\alpha_a x_{i_a})^{r_a} \\
& = & \epsilon^\alpha_j \ \sum_{ \sum r_k =  (a-1)(l-1)+r+1, r_j > 0 } (\epsilon^\alpha_1 x_{i_1})^{r_1} \dots  (\epsilon^\alpha_j x_{i_j})^{r_j} \dots (\epsilon^\alpha_a x_{i_a})^{r_a}. \\
\end{eqnarray*}
In fact, by the pigeonhole principle the condition $r_j>0$ in the last summation is redundant: if $r_j=0$, there must be some other $r_i$ greater than or equal to $l$. So this condition can be disregarded, leaving the summand symmetric. 

Now, since our choice of $j$ was arbitrary, we also have
$$ x_{\sigma_A(i_j)} h^\alpha_r(A) = \epsilon^\alpha_{\sigma_{I_a}(j)} \sum_{r_1+\dots+r_a} (\epsilon^\alpha_1 x_{i_1})^{r_1} \dots (\epsilon^\alpha_a x_{i_a})^{r_a}. $$
The result follows, since $\epsilon^\alpha_{\sigma_{I_a}(j)} = (-1)^{\alpha_j} \epsilon^\alpha_j$. \end{proof}

\begin{Lemma}\label{polymult} Let $A = \{ i_1,\dots,i_{a-1},k \}, B = \{k,j_2,\dots,j_b \} \subset I_d$ with $A \cap B = \{k\}$. Let $\alpha \in \eZ{a}$ and $\beta \in \eZ{b}$, and $r,s \geq 0$. Write $A \cup B$ for the ordered set $\{ i_1,\dots,i_{a-1},k,j_2,\dots,j_{b-1} \}$.
Then we have
\begin{equation}\label{xoneeq} h^\alpha_r(A) h^\beta_s(B) = (-1)^{\alpha_a((b-1)(l-1)+s)} h^\gamma_{r+s}(A \cup B), \end{equation}
where $\gamma = (\alpha_1,\dots,\alpha_{a-1},\beta_1,\dots,\beta_{b-1},\beta_b + \alpha_a ) \in \eZ{a+b-1}$. \end{Lemma}

\begin{proof}
The left-hand side of (\ref{xoneeq}) equals
\begin{equation}\label{arghleft}  \sum (\epsilon^\alpha_1 x_{i_1})^{r_1}\dots (\epsilon^\alpha_{a-1} x_{i_{a-1}})^{r_{a-1}} (\epsilon^\alpha_a x_k)^{r_k} (\epsilon^\beta_1 x_k)^{s_k} (\epsilon^\beta_2 x_{j_2})^{s_2} \dots (\epsilon^\beta_b x_{j_b})^{s_b} \end{equation}
\begin{eqnarray*}
\text{where} \hspace{.5in} \sum_{j=1}^{a-1} r_j + r_k & = & (a-1)(l-1) + r  \hspace{.2in} \text{and} \\
s_k + \sum_{j=2}^{b} s_j  & = & (b-1)(l-1) + s  .  \end{eqnarray*}
Now, observe that
$$ \epsilon^{\gamma}_i = \left\{ \begin{array}{ccl} 
       \epsilon^\alpha_i  & & i = 1,\dots,a \\
       &&\\
       \epsilon^\alpha_a \epsilon^\beta_{i-a} & & i = a+1,\dots,a+b-1 . \end{array} \right. $$  
Using this, we can write
\begin{equation}\label{arghright} h^\gamma_{r+s}(A \cup B) = \sum_{\begin{array}{c} r_1 , \dots , r_{a-1} \\ s_1,\dots,s_{b-1} \end{array}} \Big( \ \prod_{j=1}^{a-1} (\epsilon^\alpha_j x_{i_j})^{r_j} \prod_{j=2}^{b} (\epsilon^\beta_j x_{i_j})^{s_j} \cdot (\epsilon^\beta_1 x_k)^t (\epsilon^\alpha_a)^{t + \sum_{j=2}^{b} s_j} \ \Big) \end{equation}
where $$\sum_{j=1}^{j=a-1} r_j + \sum_{j=2}^{j=b} s_j + t = (a+b-2)(l-1)+ r + s . $$ 
Modulo a sign, every term in the summation (\ref{arghright}) arises as a term from (\ref{arghleft}), because the $r_k + s_k$ arising from (\ref{arghleft}) all satisfy the condition on $t$. On the other hand, every term in (\ref{arghleft}) must occur in (\ref{arghright}) by the pigeonhole principle: if the powers of $x_{i_1},\dots,x_{a-1}$ in a term from (\ref{arghright}) didn't sum to at least $(a-2)(l-1)+r$, then some subsequent $x_j$ would have to be raised to at least the power $l$. 

To establish the result, it remains to show that the difference in sign between the terms in (\ref{arghleft}) and (\ref{arghright}) is as claimed. (Most urgently, we need to show that it is not dependent on the choice of term.) Fixing a choice of $r_1,\dots,r_{a-1},s_2,\dots,s_b$, the ratio of the corresponding terms in (\ref{arghright}) and (\ref{arghleft}) is
$$ \big( (\epsilon^\beta_1)^t (\epsilon^\alpha_a)^{t + \sum_{j=2}^{j=b} s_j} \big) \Big/ \big( (\epsilon^\alpha_a)^{r_k} (\epsilon^\beta_1)^{s_k} \big). $$
Observing that $\epsilon^\beta_1 = 1$, we obtain
\begin{eqnarray*} \text{ratio of signs} & = & (\epsilon^\alpha_a)^{(a-1)(l-1) + r - \sum_{j=1}^{j=a-1}r_j + (b-1)(l-1) + s}(\epsilon^\alpha_a)^{r_k}  \\
& = & (\epsilon^\alpha_a)^{(b-1)(l-1)+s} \\
& = & (-1)^{\alpha_a((b-1)(l-1)+s)} . \end{eqnarray*}   \end{proof}

\begin{Lemma}\label{polymult2}
Let $l$ be odd. Let 
$$ A = \{i_1,\dots,i_{a-2},k,m \}, B= \{ k, m, j_3,\dots, j_b \}$$ 
be subsets of $I_d$ with $A \cap B = \{k,m\}$. Let $\alpha \in \eZ{a}$ and $\beta \in \eZ{b}$. Write $F$ for the polynomial $\prod_{i \in A\cup B} x_i^{l-1}$. Then
$$ h_0^\alpha(A) h_0^\beta(B) = \left\{ \begin{array}{ccl} 
       l \ F  & & \text{if  } (-1)^{\alpha_a + \beta_1} = 1, \\
       &&\\
       F & & \text{otherwise.} \end{array} \right. $$  
\end{Lemma}

\begin{proof}
By the pigeonhole principle, $h_0^\alpha(A) h_0^\beta(B)$ equals
$$ \sum_{r=0}^{l-1} x_{i_1}^{l-1}\dots x_{i_{a-2}}^{l-1} x_{j_3}^{l-1} x_{j_b}^{l-1} (\epsilon^\alpha_{a-1} x_k)^r (\epsilon^\alpha_a x_m)^{l-1-r} (\epsilon^\beta_1 x_k)^{l-1-r} (\epsilon^\beta_2 x_m)^r .$$
This is a sum of $l$ terms, each of which equals $\pm F$. Since $l-1-r$ has the same parity as $r$, and $\epsilon^\beta_1 = 1$, this sum is alternating if and only if 
$$ \epsilon^\alpha_{a-1} \epsilon^\alpha_a \epsilon^\beta_2 = -1 .$$
Observing that the left-hand side equals $(-1)^{\alpha_a + \beta_1}$ gives the result. \end{proof}

\begin{Remark}{\rm 
Retaining the notation of Lemma~\ref{polymult}:
\begin{enumerate}
\item[(i)] If $|A \cap B| = 2$ then the product $h_r^\alpha(A) h_s^\beta(B)$ is either zero or a scalar multiple of the maximal degree term $\prod_{i \in A \cup B} x_i^{l-1}$.  
\item[(ii)] If $|A \cap B | > 2$ then $h_r^\alpha(A) h_s^\beta(B) = 0$. \end{enumerate} 
}\end{Remark}

\section{X-cycles}

In this section we study the following subalgebra of $\gS_d$:

\begin{Definition}{\rm
Let $\C_X$ be the centralizer of the polynomial subalgebra of $\gS_d$. That is,
$$ \C_X = \Cen_{\gS_d}(x_1,\dots,x_d) .$$
 }\end{Definition}

\subsection{X-cycles}

\begin{Definition}\label{xcycle}{\rm
For an ordered subset $A$ of $I_d$, $r \geq 0$, and $\alpha \in \eZ{a}$, define the \emph{X-cycle} $A^{(r,\alpha)}$ by 
$$ A^{(r,\alpha)} := h_r^\alpha (A) \ \sigma_A \ c_\alpha(A) \in \gS_d .$$}\end{Definition}

\begin{Lemma} Take $A = \{i_1,\dots,i_a \} \subseteq I_d$, $\alpha \in \Z^a$, and $r \geq 0$. Then $A^{(r,\alpha)} \in \C_X$. \end{Lemma}

\begin{proof}
It is sufficient to check that $A^{(r,\alpha)}$ commutes with $x_{i_j}$ for $j \in I_a$. Using Lemma~\ref{commhx}, we have
\begin{eqnarray*}
x_{i_j} h^\alpha_r(A) \sigma_A c_\alpha(A) & = & (-1)^{\alpha_j} h^\alpha_r(A) x_{\sigma_A(i_j)} \sigma_A c_\alpha(A) \\
& = & (-1)^{\alpha_j} h^\alpha_r(A) \sigma_A x_{i_j} c_\alpha(A) \\
& = & (-1)^{\alpha_j} h^\alpha_r(A) \sigma_A (-1)^{\alpha_j} c_\alpha(A) x_{i_j} \\
& = & A^{(r,\alpha)} x_{i_j} . \end{eqnarray*}\end{proof}

For $A = \{ i_1,\dots,i_a \} \subseteq I_d$ and $\alpha \in \eZ{a}$, write $\Acy$ and $\acy$ for the elements obtained by shifting $A$ and $\alpha$ one place to the right. That is:
$$ \Acy = \{ i_a, i_1, \dots, i_{a-1} \}, \hspace{.5in} \acy = (\alpha_a,\alpha_1,\dots,\alpha_{a-1}) .$$

\begin{Lemma}\label{xshift}
For $A$ and $\alpha$ as above, and $r \geq 0$, we have
$$ \Acy^{(r,\acy)} = (-1)^{\alpha_a( (a-1)(l-1) + r + 1)} A^{(r,\alpha)} .$$
\end{Lemma}

\begin{proof}
First observe that
$$ \epsilon^{\acy}_i = \left\{ \begin{array}{ccccl} 
       \prod_{j<i-1} (-1)^{\alpha_j} \cdot \epsilon^\alpha_a & = & \epsilon^\alpha_{i-1} \epsilon^\alpha_a & & 1 < i \leq a \\
       &&&&\\
       && 1 && i=1 \end{array} \right. $$
where we use the fact that $\epsilon^\alpha_a = (-1)^{\alpha_a}$ since $\alpha$ is even. 

Now, 
$$h_r^\alpha (A) = \sum_{r_1,\dots,r_a} (\epsilon^\alpha_1 x_{i_1})^{r_1} \dots (\epsilon^\alpha_a x_{i_a})^{r_a} $$
and
$$ h_r^{\acy} (\Acy) =  \sum_{s_1,\dots,s_a} (\epsilon^\acy_1 x_{i_a})^{s_1} (\epsilon^\acy_2 x_{i_1} )^{s_2} \dots (\epsilon^\acy_a x_{i_{a-1}})^{s_a} $$
where $\sum_{j \in I_a} r_j = \sum_{j \in I_a} s_j = (a-1)(l-1) + r$. 
The second expression can be rewritten as
\begin{eqnarray*} h_r^\acy (\Acy) & = & \sum_{r_1,\dots,r_a} (\epsilon^\alpha_1 x_{i_a})^{r_a}(\epsilon^\alpha_a \epsilon^\alpha_1 x_{i_1})^{r_1} \dots (\epsilon^\alpha_a \epsilon^\alpha_{a-1} x_{i_{a-1}})^{r_{a-1}} \\
& = & \sum_{r_1,\dots,r_a} (\epsilon^\alpha_a)^{r_1+\dots+r_{a-1}} x_{i_a}^{r_a}(\epsilon^\alpha_1 x_{i_1})^{r_1} \dots (\epsilon^\alpha_{a-1} x_{i_{a-1}})^{r_{a-1}} \\
& = & \sum_{r_1,\dots,r_a} (\epsilon^\alpha_a)^{(a-1)(l-1)+r} (\epsilon^\alpha_a)^{r_a} x_{i_a}^{r_a}(\epsilon^\alpha_1 x_{i_1})^{r_1} \dots (\epsilon^\alpha_{a-1} x_{i_{a-1}})^{r_{a-1}} \\
& = & (\epsilon^\alpha_a)^{(a-1)(l-1)+r} h_r^\alpha (A) \\
& = & (-1)^{\alpha_a((a-1)(l-1)+r)} h_r^\alpha (A). \\
\end{eqnarray*}
Finally, observe that $\sigma_{\Acy} = \sigma_A$ and that $c_{\acy}(\Acy) = (-1)^{\alpha_a} c_\alpha(A)$. Putting this all together we obtain
\begin{eqnarray*}
\Acy^{(r,\acy)} & = & h^\acy_r(\Acy) \sigma_\Acy c_\acy(\Acy) \\
& = & (-1)^{\alpha_a((a-1)(l-1)+r)} h_r^\alpha (A) \sigma_A (-1)^{\alpha_a} c_\alpha(A) \\
& = & (-1)^{\alpha_a( (a-1)(l-1) + r + 1)} A^{(r,\alpha)}. \end{eqnarray*} \end{proof}

\begin{Remark}\label{xshiftrem}{\rm
Lemma~\ref{xshift} is the reason that we have to be careful to define X-cycles in terms of ordered indexing sets $A$, rather than the elements $\sigma_A$ of the symmetric group that they determine. This is quite unnatural, but it will cease to be a concern in the next section once we prove Theorem~\ref{cxshift}.  }\end{Remark}

\begin{Proposition}\label{xcymult}{\rm  
\begin{enumerate}
\item[(i)] X-cycles associated to disjoint indexing sets commute with one another.
\item[(ii)] Let $A=\{i_1,\dots,i_{a-1},k\},B=\{k,j_2,\dots,j_b\} \subset I_d$ with $A \cap B =\{k\}$, and let $r,s \geq 0$. Let $\alpha \in \eZ{a}$ and $\beta \in \eZ{b}$. Then
$$ A^{(r,\alpha)} B^{(s,\beta)} = (-1)^{\alpha_a((b-1)(l-1)+s+\beta_b)} (A \cup B)^{(r+s,\gamma)} $$
where $A \cup B = \{i_1,\dots,i_{a-1},k,j_2,\dots,j_b \}$ and
$$ \gamma = (\alpha_1,\dots,\alpha_{a-1},\beta_1,\dots,\beta_{b-1},\alpha_a+\beta_b) \in \eZ{a+b-1} .$$

\end{enumerate}

}\end{Proposition}

\begin{proof}
\begin{enumerate}
\item[(i)] Obvious.
\item[(ii)] 
\begin{eqnarray*}
A^{(r,\alpha)}B^{(s,\beta)} & = & h^\alpha_r(A) \sigma_A c_\alpha(A) h^\beta_s(B) \sigma_B c_\beta(B) \\
& = &  h^\beta_s(B) h^\alpha_r(A) \sigma_A c_\alpha(A) \sigma_B c_\beta(B) \hspace{.5in}  \\
& = & h^\alpha_r(A) h^\beta_s(B) \sigma_A \sigma_B c_\alpha(\{i_1,\dots,i_{a-1},j_b\}) c_\beta(B) \\
& = & (-1)^{\alpha_a((b-1)(l-1)+s)} h^\gamma_{r+s}(A \cup B) \sigma_{A \cup B} \cdot \\
& & \hspace{1.5in} c_\alpha(\{i_1,\dots,i_{a-1},j_b\}) c_\beta(B)  \\ 
& = & (-1)^{\alpha_a((b-1)(l-1)+s)+\alpha_a \beta_b} h^\gamma_{r+s}(A \cup B) \sigma_{A \cup B} c_\gamma (A \cup B) \end{eqnarray*}
from which the result follows. \end{enumerate} \end{proof}

\subsection{Other elements of $\C_X$}

\begin{Lemma}\label{evenx}
Take $A = \{ i_1,\dots,i_a \} \subseteq I_d$ and $\alpha \in \Z^a$. Suppose 
$$ z = f \sigma_A c_\alpha(A) \in \C_X$$ 
for some polynomial $f \in R_l[A]$. Then either $\alpha \in \eZ{a}$, or $f$ is of maximal degree with respect to $A$.
\end{Lemma}

\begin{proof}
By multiplying $z$ on either side by $x_{i_j}$, for $j \in I_a$, we see that
$$ x_{i_j} f = (-1)^{\alpha_j} x_{\sigma_A(i_j)} f . $$
This implies
$$ x_{i_1} f = \prod_{i \in I_a} (-1)^{\alpha_i} x_{i_1} f , $$
so either $\prod_{i \in I_a} (-1)^{\alpha_i} = 1$, whereupon $\alpha \in \eZ{a}$, or else $x_{i_1} f = 0$. In the latter scenario, we must then have $x_{i_j} f = 0$ for all $j \in I_a$, which gives the required condition on $f$.\end{proof}

\begin{Lemma}\label{evenx2} Let $A = \{i_1,\dots,i_a\} \subseteq I_d$ and $\alpha \in \eZ{a}$. Suppose $z = f \sigma_A c_\alpha(A) \in \C_X$, for some homogeneous $f \in R_l[A]$ with $\deg(f) < (a-1)(l-1)$. Then $z$ is a scalar multiple of an X-cycle.   \end{Lemma}

\begin{proof}
By assumption, there exists a $j$ with $x_{i_j} f \neq 0$. Multiplying $z$ on both sides by $x_{i_j}$, we see that $f$ must satisfy the relation (\ref{commhxform}) for this choice of $j$, whereupon $x_{\sigma_A(i_j)} f \neq 0$. The result follows from repeated application of this argument, and the observation that (\ref{commhxform}) characterizes the polynomial $h_r^\alpha(A)$ up to a scalar. \end{proof}

\begin{Proposition} Let $z = f \sigma_{A_1} \dots \sigma_{A_m} c \in \C_{x}$, for some pairwise disjoint $A_1,$ $\dots,A_m$ in $I_d$, $f$ homogeneous in $R_l[A_1 \cup \dots \cup A_m]$ and $c \in C[I_d]$. Then $z$ can be written as a product of disjoint X-cycles and maximal degree cycles.  \end{Proposition}

\begin{proof}
Proceed by induction on $m$, the case $m=1$ following from Lemma~\ref{evenx} and Lemma~\ref{evenx2}. It is clear from the relations that we may assume $c \in C[A_1 \cup \dots \cup A_m]$.  

Write $f = \sum_i f_i g_i$ for some homogeneous polynomials $f_i \in R_l[A_1 \cup \dots A_{m-1}]$ and $g_i \in R_l[A_m]$. The element $c$ can be written uniquely in the form 
$$c_\alpha(A_1 \cup \dots \cup A_{m-1}) c_\beta(A_m).$$ 
Now, consider the coefficient of each $g_i \sigma_{A_m} c_\beta(A_m)$ in $z$: by homogeneity of the relations, each such coefficient must commute with $R_l[A_1 \cup \dots \cup A_{m-1}]$, and so by the inductive hypothesis it is of the required form. Moreover, each $g_i \sigma_{A_m} c_\beta(A_m)$ must commute with $R_l[A_d]$ so, depending on the parity of $\beta$, apply Lemma~\ref{evenx} and/or Lemma~\ref{evenx2} to obtain the result. \end{proof}

\section{CX-cycles}

\begin{Definition}{\rm 
Let $\C_{CX}$ be the even centralizer in $\C_X$ of the Clifford subalgebra of $\gS_d$.  
}\end{Definition}

\subsection{CX-cycles} In order to refine the notion of X-cycle and obtain elements of $\C_{CX}$, we need to introduce another sign function on $\eZ{a}$. 

\begin{Definition}{\rm
For $\alpha \in \eZ{a}$, define
$$ \tau_\alpha := (-1)^{\frac{1}{2} |\alpha| + \sum_{i \in I_a} i \alpha_i } .$$ }\end{Definition}

Recall the notation $\alpha^{(j)}$ from Section~\ref{notation}. The following is another (and less ad hoc) way of characterizing the $\tau_\alpha$:

\begin{Lemma}\label{taulem} Let $\alpha \in \eZ{a}$ and $1 < j \in I_a$.
\begin{enumerate}
\item[(i)] $\tau_{\alpha^{(j)}} = (-1)^{\alpha_j + \alpha_{j-1}} \tau_\alpha$.
\item[(ii)] $\tau_{\alpha^{(1)}} = (-1)^{\alpha_1 + \alpha_a + a} \tau_\alpha$. 
\end{enumerate}
\end{Lemma}

\begin{proof}
\begin{enumerate}
\item[(i)] Observe that
$$ |\alpha^{(j)}|  = \left\{ \begin{array}{ccl} 
       |\alpha| \pm 2 & & \text{if} \ \alpha_j = \alpha_{j-1}, \\
       &&\\
       |\alpha| && \text{otherwise.} \end{array} \right. $$
Also
$$ \sum_{i \in I_a} i (\alpha^{(j)})_i \equiv \sum_{i \in I_a} i \alpha_i + j + (j+1) \text{  mod } 2 $$
and so
$$ \tau_{\alpha^{(j)}}  = \left\{ \begin{array}{ccl} 
       \tau_\alpha & & \text{if} \ \alpha_j = \alpha_{j-1}, \\
       &&\\
       - \tau_\alpha && \text{otherwise.} \end{array} \right. $$
The result follows, since $(-1)^{\alpha_j + \alpha_{j-1}}$ is negative precisely when $\alpha_j \neq \alpha_{j-1}$. 

\item[(ii)] Similarly to part (i), we have
$$ |\alpha^{(1)}|  = \left\{ \begin{array}{ccl} 
       |\alpha| \pm 2 & & \text{if} \ \alpha_1 = \alpha_a, \\
       &&\\
       |\alpha| && \text{otherwise,} \end{array} \right. $$
but then we have
$$ \sum_{i \in I_a} i (\alpha^{(1)})_i \equiv \sum_{i \in I_a} i \alpha_i + 1 + a  \text{  mod } 2$$
and so
$$ \tau_{\alpha^{(1)}}  = \left\{ \begin{array}{ccl} 
       (-1)^a \tau_\alpha & & \text{if} \ \alpha_1 = \alpha_a, \\
       &&\\
       (-1)^{a+1} \tau_\alpha && \text{otherwise.} \end{array} \right. $$ \end{enumerate} \end{proof}

\begin{Definition}{\rm 
\begin{enumerate}
\item[(i)] For $A = \{i_1,\dots,i_a \} \subseteq I_d$ and $r \geq 0$, define the element $A^{(r)}$ as follows:
\begin{equation}\label{cxcycleform} A^{(r)} := \sum_{\alpha \in \eZ{a}} \tau_\alpha A^{(r,\alpha)} \in \C_X  \end{equation}
\item[(ii)] Call such an element a \emph{CX-cycle} if $(a-1)l+r$ is even.
\end{enumerate}
}\end{Definition}

\begin{Proposition}{
The element $A^{(r)}$ defined above lies in $\C_{CX}$ if and only if it is a CX-cycle.
}\end{Proposition}

\begin{proof}
It is obviously sufficient to check that $A^{(r)}$ commutes with $c_{i_j}$ for $j \in I_a$. We treat the cases $j=1$ and $j > 1$ separately.

First, suppose $j >1$. Then
\begin{eqnarray*} c_{i_j} A^{(r)} & = & \sum_{\alpha \in \eZ{a}} \tau_\alpha c_{i_j} h^\alpha_r(A) \sigma_A c_\alpha(A) \\
& = & \sum_{\alpha \in \eZ{a}} \tau_{\alpha} h^{\alpha^{(j)}}_r(A) c_{i_j} \sigma_A c_\alpha(A) \hspace{.8in} \text{by Lemma~\ref{commhc}~(i)} \\
& = & \sum_{\alpha \in \eZ{a}} \tau_{\alpha} h^{\alpha^{(j)}}_r(A) \sigma_A c_{i_{j-1}} c_\alpha(A) \\
& = & \sum_{\alpha \in \eZ{a}} \epsilon^\alpha_{j-1} \tau_{\alpha} h^{\alpha^{(j)}}_r(A) \sigma_A c_{\alpha+ 1_{j-1}}(A) \hspace{.4in} \text{by Remark~\ref{signprop}~(ii).} \end{eqnarray*}
On the other hand,
\begin{eqnarray*}
 A^{(r)} c_{i_j} & = & \sum_{\alpha \in \eZ{a}} \tau_\alpha h^\alpha_r(A) \sigma_A c_\alpha(A) c_{i_j} \\
& = & \sum_{\alpha \in \eZ{a}} \epsilon^\alpha_{j+1} \tau_\alpha  h^\alpha_r(A) \sigma_A c_{\alpha+1_j} (A) \hspace{.8in} \text{by Remark~\ref{signprop}~(ii)} \\
& = & \sum_{\alpha \in \eZ{a}} \epsilon^{\alpha^{(j)}}_{j+1} \tau_{\alpha^{(j)}} h^{\alpha^{(j)}}_r(A) \sigma_A c_{\alpha + 1_{j-1}} (A) \hspace{.4in}  \\
& = & \sum_{\alpha \in \eZ{a}} \epsilon^{\alpha^{(j)}}_{j+1} (-1)^{\alpha_j+\alpha_{j-1}} \tau_\alpha h^{\alpha^{(j)}}_r(A) \sigma_A c_{\alpha+1_{j-1}} (A) \\
& = & c_{i_j} A^{(r)}  \hspace{2in} \text{because } \epsilon^{\alpha^{(j)}}_{j+1} = \epsilon^\alpha_{j+1}.  \end{eqnarray*}

Now suppose $j=1$. We have
\begin{eqnarray*} 
c_{i_1} A^{(r)} & = & \sum_{\alpha \in \eZ{a}} \tau_\alpha c_{i_1} h^\alpha_r(A) \sigma_A c_\alpha(A) \\
& = & (-1)^{(a-1)(l-1)+r} \sum_{\alpha \in \eZ{a}} \tau_\alpha  h^{\alpha^{(1)}}_r(A) c_{i_1} \sigma_A c_\alpha(A) \\
& = & (-1)^{(a-1)(l-1)+r} \sum_{\alpha \in \eZ{a}} \tau_\alpha h^{\alpha^{(1)}}_r(A) \sigma_A c_{i_{a}} c_\alpha(A) \\
& = & (-1)^{(a-1)(l-1)+r} \sum_{\alpha \in \eZ{a}} (-1)^{\alpha_a} \tau_\alpha h^{\alpha^{(1)}}_r(A) \sigma_A c_{\alpha+1_a}(A), \\
\end{eqnarray*}
whereas
\begin{eqnarray*}
A^{(r)} c_{i_1} & =  & \sum_{\alpha \in \eZ{a}} \tau_\alpha h^\alpha_r(A) \sigma_A c_\alpha(A) c_{i_1} \\
& = & \sum_{\alpha \in \eZ{a}} (-1)^{\alpha_1} \tau_\alpha h^\alpha_r(A) \sigma_A c_{\alpha+1_1}(A) \\
& = & \sum_{\alpha \in \eZ{a}} (-1)^{(\alpha^{(1)})_1} \tau_{\alpha^{(1)}} h^{\alpha^{(1)}}_r(A) \sigma_A c_{\alpha+1_a}(A) \\
& = & \sum_{\alpha \in \eZ{a}} (-1)^{\alpha_1+1} (-1)^{\alpha_1+\alpha_a+a} \tau_\alpha h^{\alpha^{(1)}}_r(A) \sigma_A c_{\alpha + 1_a} (A) . \\ \end{eqnarray*}
Comparing the signs, we see that $A^{(r)}$ commutes with $c_{i_1}$ if and only if 
$$ (a-1)(l-1)+r+\alpha_a +\alpha_a+a+1 $$ 
is even. This is equivalent to the condition that $(a-1)l+r$ be even, which establishes the result. \end{proof}

\begin{Theorem}\label{cxshift}
Let $A^{(r)}$ be a CX-cycle. Let $\Acy$ be defined as in Lemma~\ref{xshift}. Then we have
$$ \Acy^{(r)} = A^{(r)} .$$   \end{Theorem}

\begin{proof}
By definition, 
$$ \Acy^{(r)} = \sum_{\alpha \in \eZ{a}} \tau_\alpha \Acy^{(r,\alpha)} = \sum_{\alpha \in \eZ{a}} \tau_{\acy} \Acy^{(r,\acy)} .$$
Since $|\acy| = |\alpha|$ for any $\alpha$, we have that 
\begin{equation}\label{taus} \tau_\acy = \left\{ \begin{array}{ccl} 
       - \tau_\alpha & & \text{if} \ \alpha_a=1 \ \text{and} \ a \ \text{is odd, and} \\
       &&\\
        \tau_\alpha && \text{otherwise.} \end{array} \right. \end{equation}
Now, combining (\ref{taus}) with Lemma~\ref{xshift}, the result follows from consideration of the following three cases:
\begin{enumerate}
\item[(i)] $\alpha_a = 0$. Then $\tau_\acy \Acy^{(r,\acy)} = \tau_\alpha A^{(r,\alpha)}$ automatically.
\item[(ii)] $\alpha_a=1$ and $a$ is odd. Then
$$ \tau_\acy \Acy^{(r,\acy)} = \tau_\alpha A^{(r,\alpha)} \ \leftrightarrow \ (a-1)(l-1)+r  \text{ is even} \ \leftrightarrow \ (a-1)l + r \text{ is even.} $$
\item[(iii)] $\alpha_a=1$ and $a$ is even. Then \begin{eqnarray*}
\tau_\acy \Acy^{(r,\acy)} = \tau_\alpha A^{(r,\alpha)} & \leftrightarrow & (a-1)(l-1)+r+1 \text{ is even} \\
& \leftrightarrow & r \text{ and } l \text{ have the same parity} \\
& \leftrightarrow & (a-1)l + r \text{ is even.} \end{eqnarray*}  \end{enumerate} \end{proof}

\begin{Remark}{\rm
In contrast to Remark~\ref{xshiftrem}, the above theorem shows that CX-cycles are ``well-defined'' in a way that X-cycles are not. That is, a CX-cycle can be unambiguously associated to a cycle in the symmetric group and a positive integer. It is clear that CX-cycles associated to disjoint cycles commute with one another, and that for $A,B \subset I_d$ and $r,s \geq 0$ the product $A^{(r)} B^{(s)}$ must be zero whenever $|A \cap B| > 2$, or when $|A \cap B| = 2$ and $r+s>0$.  }\end{Remark}

\begin{Proposition}\label{cxmult1} Let $A = \{ i_1,\dots, i_{a-1}, k \}, B = \{ k, j_2,\dots,j_b \} \subseteq I_d$ with $A \cap B = \{k \}$, and let $r,s \geq 0$ be such that $B^{(s)}$ is a CX-cycle. As in Lemma~\ref{polymult} and Proposition~\ref{xcymult}, write $A \cup B$ for the ordered index set 
$$ \{i_1,\dots,i_{a-1},k,j_1,\dots,j_{b-1} \}.$$ 
Then
$$ A^{(r)} B^{(s)} = (A \cup B)^{(r+s)} .$$   \end{Proposition}

\begin{proof}
As in Proposition~(\ref{xcymult})~(ii), for given $\alpha \in eZ{a}$ and $\beta \in \eZ{b}$, set
$$ \gamma =(\alpha_1,\dots,\alpha_{a-1},\beta_1,\dots,\beta_{b-1},\alpha_a+\beta_b) \in \eZ{a+b-1} .$$
Observe that each element of $\eZ{a+b-1}$ is determined by a unique pair of elements of $\eZ{a}$ and $\eZ{b}$ in this way. Now,
\begin{eqnarray*}
A^{(r)} B^{(s)} & = & \sum_{\alpha,\beta} \tau_\alpha \tau_\beta A^{(r,\alpha)} B^{(s,\beta)} \\
 & = & \sum_{\alpha, \beta} (-1)^{\alpha_a((b-1)(l-1)+s+\beta_b)} \tau_\alpha \tau_\beta (A \cup B)^{(r+s,\gamma)}.  \end{eqnarray*}
Now we need to show that this coefficient equals $\tau_\gamma$. Observe that
$$ |\gamma| = \left\{ \begin{array}{ccl} 
       |\alpha| + |\beta| - 2  & & \text{if } \alpha_a = \beta_b = 1, \\
       &&\\
       |\alpha| + |\beta| & & \text{otherwise} . \end{array} \right. $$  
Also,
\begin{eqnarray*}
\sum_{i \in I_{a+b-1}} i \gamma_i & = & \sum_{i \in I_{a-1}} i \alpha_i + \sum_{j \in I_{b-1}} (a-1+j)\beta_j + (a+b-1)(\alpha_a+\beta_b) \\
& = & \sum_{i \in I_{a}} i \alpha_i + \sum_{j \in I_{b}} j\beta_j + (b-1)\alpha_a + (a-1) \sum_{j \in I_b} \beta_j .\end{eqnarray*}
Since the last term is even, this gives us that
$$ \tau_\gamma = \tau_\alpha \tau_\beta (-1)^{\alpha_a(\beta_b + b -1 )} .$$
Comparing this to our initial calculation, we need $b-1$ and $(b-1)(l-1)+s$ to have the same parity. However, this is equivalent to the requirement that $(b-1)l+s$ be even, so we're done.  \end{proof}

\begin{Example}\label{factor}{\rm
The requirement in Proposition~(\ref{cxmult1}) that one factor be a CX-cycle is strict. For instance, the factorization
$$ (1 \ 2 \ 3)^{(0)} = (1 \ 2)^{(0)} (2 \ 3)^{(0)} $$
is valid if $l$ is even, but not if $l$ is odd. 
}\end{Example}

\begin{Example}{\rm $ (1 \ 2)^{(0)} (1 \ 2)^{(0)} = 0$, regardless of the parity of $l$. }\end{Example}

\begin{Example}\label{simpleoverlap2}{\rm 
\begin{enumerate}
\item[(i)] As a consequence of the previous two examples, we have
$$ (1 \ 2 \ 3)^{(0)} (2 \ 3 \ 4)^{(0)} = (1 \ 2)^{(0)} (2 \ 3)^{(0)} (2 \ 3)^{(0)} (3 \ 4)^{(0)} = 0 $$
if $l$ is even. However, if $l$ is odd, one can use Lemma~\ref{polymult2} to show that
$$ (1 \ 2 \ 3)^{(0)} (2 \ 3 \ 4)^{(0)} = 2 x_1^{(l-1)} x_2^{(l-1)} x_3^{(l-1)} x_4^{(l-1)} s_1 s_3 ( c_1c_4 + c_2c_3 - c_2 c_4 - c_1 c_3   ) .$$
\item[(ii)] Regardless of the parity of $l$, we have
$$ (1 \ 2 \ 3)^{(0)} (3 \ 2 \ 4)^{(0)} = 0. $$
\end{enumerate}
}\end{Example}

\begin{Remark}{\rm
Collectively Theorem~\ref{cxshift}, Proposition~\ref{cxmult1}, and Example~\ref{simpleoverlap2} provide sufficient information to calculate $A^{(0)} B^{(0)}$ whenever $|A \cap B| = 2$. (This is only complicated if $l$ is odd.)
}\end{Remark}

\begin{Lemma}\label{cxsuff} Let $A = \{i_1,\dots,i_a \} \subseteq I_d$, and 
$$ z = \sum_{\alpha \in \eZ{a}} f_\alpha \sigma_A c_\alpha (A) \in \C_{CX} $$
for some homogeneous $f_\alpha \in R_l[A]$. If the $f_\alpha$ are not of maximal degree with respect to $A$, then $z$ is proportional to a CX-cycle.   \end{Lemma}

\begin{proof}
We may write $z = \sum_{\alpha \in \eZ{a}} \theta_\alpha A^{(r,\alpha)}$ for some $0 \leq r \leq l-1$. Multiplying this expression on the left and on the right by $c_i$, for $i \in A$, and using Lemma~\ref{commhc}, we obtain
$$ \theta_{\alpha^{(i)}} = (-1)^{\alpha_i + \alpha_{i-1}} \theta_\alpha $$
for each $1 < i \in I_a$, as well as a similar expression in the case $i=1$. However, this characterizes the coefficients $\tau_\alpha$ up to a scalar, so the result follows. \end{proof}

\subsection{Odd skew elements}

\begin{Example}{\rm If $l$ is odd, we have that 
$$ (1 \ 2 \ 3)^{(0)} (2 \ 3 \ 4)^{(0)} = 2 x_1^{(l-1)} x_2^{(l-1)} x_3^{(l-1)} x_4^{(l-1)} s_1 s_3 ( c_1c_4 + c_2c_3 - c_2 c_4 - c_1 c_3   ) .$$
The right-hand side cannot be expressed as a product of disjoint CX-cycles, but instead factors as follows:
$$ - 2 ( x_1^{(l-1)} x_2^{(l-1)} s_1 (c_1 - c_2) ) \cdot ( x_3^{(l-1)} x_4^{(l-1)} s_3 (c_3 - c_4) ).$$ 
The two terms appearing in this factorization are maximal degree odd elements of $\C_X$, each of which skew-commutes with the Clifford part of $\gS_d$. 
 }\end{Example}

\begin{Definition}{\rm Say an element of the form $f \sigma_A c \in \C_X$, for $A = \{i_1,\dots,i_a\} \subseteq I_d$, $f \in R_l[A]$, and $c \in C[A]$, is an \emph{odd skew cycle} if it is odd and skew-commutes with all the Clifford generators of $\gS_d$. }\end{Definition}

\begin{Proposition}\label{oddla} For $A= \{i_1,\dots,i_a\} \subseteq I_d$, $f \in R_l[A]$, and $c \in C[A]$, let $z = f \sigma_A c$ be an odd skew cycle. Then it is of maximal degree with respect to $A$. Moreover, one of $l$ and $a$ must be even. \end{Proposition}

\begin{proof}
Since $z$ is odd, we can write it as $\sum_{\alpha \in \oZ{a}} f_\alpha \sigma_A c_\alpha(A)$. By Lemma~\ref{evenx}, the $f_\alpha$ must all be of maximal degree with respect to $A$, so we can in fact write
$$ z = F \sigma_A \big( \sum_{\alpha \in \oZ{a}} \theta_\alpha c_\alpha(A) \big) $$
for some coefficients $\theta_\alpha$, where $F = \prod_{j=1}^a x_{i_j}^{l-1}$.
Let $1<j \in I_a$. Multiplying on the left by $c_{i_j}$, we obtain
$$ (-1)^{l-1} \sum_{\alpha \in \oZ{a}} \big( \theta_\alpha \prod_{k=1}^{j-2} (-1)^{\alpha_k} c_{\alpha+(j-1)}(A) \big) $$
whereas multiplying on the right we obtain
$$ \sum_{\alpha \in \oZ{a}} \big( \theta_\alpha \prod_{k=j+1}^a (-1)^{\alpha_k} c_{\alpha+j}(A) \big) .$$
Observing that if $\alpha \in \oZ{a}$ we have
$$ \prod_{k=1}^{j-2} (-1)^{\alpha_k} \cdot \prod_{k=j+1}^a (-1)^{\alpha^{(j)}_k} = - (-1)^{\alpha_j + \alpha_{j-1}}, $$ 
and comparing coefficients yields
\begin{equation}\label{thetaform} \theta_{\alpha^{(j)}} = (-1)^{l-1} (-1)^{\alpha_j + \alpha_{j-1}} \theta_\alpha . \end{equation}
One can similarly obtain
$$ \theta_{\alpha^{(1)}} = (-1)^{l-1} (-1)^{\alpha_1 + \alpha_{a}} \theta_\alpha .$$
Now, using the fact that $\alpha = (\dots(\alpha^{(a)})^{(a-1)}\dots)^{(1)}$, we have
\begin{eqnarray*} \theta_\alpha & = & (-1)^{a(l-1)} (-1)^{\alpha_a+\alpha_{a-1}} (-1)^{\sum_{j=2}^{a-1} (\alpha_j + \alpha_{j-1} + 1)} (-1)^{\alpha_1+\alpha_a + 2} \theta_\alpha \\
& = & (-1)^{al} \theta_\alpha,  \end{eqnarray*}
from which the claim follows.   \end{proof}

\subsection{Basis for $\C_{CX}$}\label{cxbasis}

\begin{Theorem}
\begin{enumerate}
\item[(i)] Let $A_1,\dots ,A_m$ be pairwise disjoint subsets of $I_d$, and let $z = F \sigma_{A_1} \dots \sigma_{A_m} c \in \C_{CX}$ be such that $c \in C_{A_1 \cup \dots \cup A_m}$, and $F \in R_l[I_d]$ is of maximal degree with respect to $A_1 \cup \dots \cup A_m$. Then $z$ can be written as a product of disjoint even CX-cycles and (an even number of) disjoint odd skew cycles. 
\item[(ii)] Let $z = \sum_{\alpha \in \eZ{a}} f_\alpha \sigma_{A_1} \dots \sigma_{A_m} c_\alpha(A_1 \cup \dots \cup A_m)$ be an element of $\C_{CX}$, where $A_1,\dots,A_m$ are pairwise disjoint subsets of $I_d$, $|A_i| = a_i$ for $i \in I_m$, and $\sum_{i \in I_m} a_i = a$. Suppose the $f_\alpha$ are not of maximal degree with respect to any of the $A_i$. Then $z$ can be written as a sum of products of disjoint CX-cycles.
\end{enumerate}
\end{Theorem}

\begin{proof} We prove the first statement by induction on $m$. For the base case, it is easy (using Lemma~\ref{taulem}) to calculate that $z$ must equal the CX-cycle $A_1^{(l-1)}$. (And thereby that $|A_1|$ must be even if $l$ is even.)

For the induction step, begin by writing 
$$c = \sum_{\alpha,\beta} \theta_{\alpha,\beta} c_\alpha(A_1 \cup \dots \cup A_{m-1}) c_\beta(A_m) $$
for appropriate $\alpha, \beta$ and some coefficients $\theta_{\alpha,\beta}$. Consider the coefficient of each $\prod_{i \in A_m}x_i^{l-1} \sigma_{A_m} c_\beta(A_m)$ in $z$. 

Firstly, if $\beta$ is even then this coefficient must itself be even, and must commute with every $c_i$ with $i \in A_1 \cup \dots \cup A_{m-1}$. So, by the inductive hypothesis, it is of the appropriate form. If $\beta$ is odd, then the coefficient is odd and skew-commutes with all these elements. If $m=2$, this means that the coefficient is a single odd skew cycle. Otherwise, repeat this argument until the inductive hypothesis can be applied. Finally, by grouping like coefficients and appealing to homogeneity of the relations, we reduce to the base case. 

The second statement is proved similarly, using Lemma~\ref{cxsuff} for the base case and Proposition~\ref{oddla} to observe that no odd elements can arise in the factorization.  \end{proof}

\section{Center of $\gS_d$}\label{penultimate}

\subsection{Standard basis}

Let $\M_d(l)$ denote the set of $l$-multipartitions of $d$.

\begin{Definition}{\rm 
Let $z = A_1^{(r_1)} \dots A_m^{(r_m)}$ be a product of disjoint CX-cycles in $\C_{CX}$. Define the \emph{cycle type} of $z$ to be the multipartition ${\bf \lambda} = (\lambda^{(1)},\dots,\lambda^{(l)}) \in \M_d(l)$, where the parts of $\lambda^{(r)}$ consist of those $|A_i|$ occurring in $z$ such that $r_i = r-1$. 
}\end{Definition}

\begin{Definition}{\rm
For $\lambda \in \M_d(l)$, define the \emph{redundancy} of $ \lambda$ to be the number of parts of $\lambda^{(1)}$ that are equal to $1$. So if $z$ is of cycle type $\lambda$, then the redundancy of $\lambda$ is the number of $1$-cycles in $z$ of degree $0$. (These could be omitted from $z$ because $(i)^{(0)} = 1$.)
}\end{Definition}

The next lemma is an immediate consequence of the definitions:

\begin{Lemma}\label{mevlem} Let $z$ be defined as above, and let ${\bf \lambda}$ be the cycle type of $z$. Then one of the following is true:
\begin{enumerate}
\item[(i)] $l$ is even. $\lambda^{(r)} = \emptyset$ whenever $r$ is even.
\item[(ii)] $l$ is odd. $\lambda^{(r)}$ has only odd parts if $r$ is odd, and only even parts if $r$ is even.  
\end{enumerate}
\end{Lemma}

Let $\eM_d(l) \subset \M_d(l)$ consist of those multipartitions which occur as the cycle type of some product of disjoint CX-cycles. 

\begin{Definition}{\rm
For ${\bf \lambda} \in \eM_d(l)$, let $z_d({\bf \lambda})$ be the sum of all CX-cycles with cycle type ${\bf \lambda}$.
}\end{Definition}

\begin{Theorem}\label{zbasis} Let $l$ be odd. Then the elements $\{ \ z_d({\bf \lambda}) \ | \ {\bf \lambda} \in \eM_d(l) \ \}$ form a basis for $Z( \gS_d )$. \end{Theorem}

The symmetric group $\Sigma_d$ acts on $\C_{CX}$ by conjugation, and $Z(\gS_d)$ consists of the fixed points of this action. We consider the orbit sums of elements of $\C_{CX}$. It is clear that $z_d(\lambda)$ is the orbit sum corresponding to any product of CX-cycles with cycle type $\lambda$. So, in view of section~\ref{cxbasis}, the theorem follows from the following lemma.

\begin{Lemma} Let $l$ be odd, and let $z = F \sigma_A c$ be an odd skew cycle, for some $A = \{i_1,\dots,i_a \} \subset I_d$, $F \in R_l[A]$, $c \in C[A]$. Then 
$$ \sum_{\tau \in \Sigma_a} (\tau \cdot z ) = 0.$$  \end{Lemma}

\begin{proof}
From Proposition~\ref{oddla}, we have that $a$ is even. Write
$$ \sum_{\tau \in \Sigma_a} (\tau \cdot z ) = \sum_{i=1}^{(a-1)!} \sum_{\tau \in \tau_i <\sigma_A>} (\tau \cdot z),$$
where $\{\tau_i\}$ is a set of coset representatives for $<\sigma_A>$ in $\Sigma_a$. We claim that the sum over each coset is zero. As in Proposition~\ref{oddla}, write 
$$ c = \sum_{\alpha \in \oZ{a}} \theta_\alpha c_\alpha(A),$$ 
and consider the element 
$$\sigma_A^{-1} c \sigma_A(A)= \sum_{\alpha \in \oZ{a}} \theta_{\sigma_A^{-1}(\alpha)} c_\alpha(A) \in C[A]. $$ 
We can obtain $\sigma_A^{-1}(\alpha)$ from $\alpha$ by applying a sequence of precisely $|\alpha|$ operations in the following manner: find $i \in I_a$ such that $\alpha_i = 1$ and $\alpha_{i-1}=0$. (This is possible since $a$ is even and $\alpha \in \oZ{a}$.) Replace $\alpha$ with $\alpha^{(i)}$, then continue sequentially. Now, applying (\ref{thetaform}), we have $\theta_{\sigma_A^{-1}(\alpha)} = - \theta_\alpha$. Consequently, 
$$ \sigma_A \cdot z = - z ,$$
so we have $\sum_{\tau \in <\sigma_A>} (\tau \cdot z) = 0$, and hence the result.  \end{proof}

\subsection{Murphy basis}

Now we construct another basis for $Z(\gS_d)$. Its indexing set, given in Definition~\ref{pev}, has a more natural description than the one in Lemma~\ref{mevlem}.

For a partition $\lambda = (\lambda_1 \geq \dots \geq \lambda_r > 0)$, define 
$$\lambda/l := (\lfloor \lambda_1/l \rfloor \geq \dots \lfloor \lambda_r/l \rfloor ).$$

\begin{Definition}\label{pev}{\rm 
\begin{enumerate}
\item[(i)] Let $\P_d(l)$ denote the set of partitions $\lambda$ satisfying 
$$r + |\lambda/l| \leq d. $$ 
\item[(ii)] Let $\eP_d(l)$ be the subset of $\P_d(l)$ consisting of partitions all of whose parts are even.
\end{enumerate} }\end{Definition}

We use the bijection $\varphi: \M_d(l) \rightarrow \P_d(l)$ from section $2$ of \cite{B}. If 
$$ {\bf \lambda} = (\lambda^{(1)},\dots,\lambda^{(l)}) \in \M_d(l), $$ 
where $\lambda^{(r)} = (\lambda^{r}_1 \geq \dots \geq \lambda^{(r)}_{m_r} > 0)$, then $\varphi({\bf \lambda})$ is the partition with parts $(\lambda^{(r)}_i - 1)l + r - 1$ for all $ r \in I_l$ and $1 \leq i \leq m_r$. 

(In other words: the $r$th entry of $\lambda$ encodes the parts of $\varphi({\bf \lambda})$ that are congruent to $r-1$ mod $l$.)

\begin{Lemma} 
\begin{enumerate}
\item[(i)] $\varphi$ is a bijection.
\item[(ii)] $\varphi(\eM_d(l)) = \eP_d(l)$. 
\end{enumerate}
\end{Lemma}

\begin{proof}
The construction of $\phi^{-1}$ is also given in \cite{B}: given $\mu \in \P_d(l)$, construct ${\bf \lambda} = (\lambda^{(1)},\dots,\lambda^{(l)})$ as follows. For all $\mu_i$ satisfying $\mu_i \equiv r-1$ mod $l$, $\lambda^{(r)}$ has parts $\{ \lfloor \mu_i / l \rfloor + 1 \}$. The second part is a direct consequence of Lemma~\ref{mevlem}. \end{proof}

\begin{Definition}{\rm
For $k = (a-1)l + r \geq 0 $, define the \emph{colored Jucys-Murphy element} 
\begin{equation}\label{yik} y_i(k) = \sum_{ \begin{array}{c}
1 \leq i_1,\dots,i_{a-1} < i \\
i_1, \dots, i_{a-1} \ \text{distinct} \end{array}} (i \ i_1 \dots i_{a-1} )^{(r)} \in \gS_d. \end{equation}
}\end{Definition}

\begin{Remark}{\rm
\begin{enumerate}
\item[(i)] The terms on the right-hand side of (\ref{yik}) are CX-cycles if and only if $k$ is even. 
\item[(ii)] If $k<l$ we have 
$$ y_i(k) = x_i^k .$$ 
\item[(iii)] Something that is to be called a ``Jucys-Murphy element'' should be expected to consist of a sum of elements looking (in this instance) like 
$$(i_1 \ \dots \ i_{a-1} \ i )^{(r)}.$$ 
However, because of Theorem~\ref{cxshift}, our definition is equivalent to this if $k$ is even. (But not if $k$ is odd.) The next lemma explains the irregularity. 
\end{enumerate}
}\end{Remark}

\begin{Lemma}\label{lcase}
For all $i \in I_d$ we have $\x_i^l \in F_{l-1} \S_d$, and $\gr_{l-1} \x_i^l = y_i(l)$.  
\end{Lemma}

\begin{proof}
Proceed by induction on $i$, the case $i=1$ being clear since $y_1(l)=0$ and the polynomial $f$ has no degree $l-1$ term. 

As a consequence of (\ref{affpowerrel}), we have 
$$ \x_{i+1}^l = s_i \x_i^l s_i + \sum_{j=0}^{l-1} \x_i^j \x_{i+1}^{l-1-j} \s_i (1 + (-1)^{j+1} \c_i \c_{i+1}) $$
in $\S_d$, and the image of the summation under $\gr_{l-1}$ can be seen to be equal to $(i \ i+1)^{(0)}$. So, by the inductive hypothesis, $\x_{i+1}^l \in F_{l-1} \S^f_d$, and
\begin{eqnarray*}
\gr_{l-1} \x_{i+1}^l & = & s_i y_i(l)  s_i + (i \ i+1)^{(0)} \\
& = & y_{i+1}(l)  . \end{eqnarray*} \end{proof}

We can regard elements of $\eP_d(l)$ simply as $d$-tuples of integers. For two such $d$-tuples, $\mu$ and $\nu$, write $\mu \sim \nu$ if they differ up to permutation of their entries. 

\begin{Definition}{\rm
For $\mu \in \eP_d(l)$, define
$$ m_d(\mu) := \sum_{(\nu_1,\dots,\nu_d) \sim \mu} y_1(\nu_1) \dots y_d(\nu_d) .$$
}\end{Definition}

\begin{Theorem}\label{mbasis} Let $l$ be odd. Then the elements $\{ \ m_d(\mu) \ | \ \mu \in \eP_d(l) \ \}$ form a basis for $Z(\gS_d)$.  \end{Theorem}

\begin{proof}
First we need to show that each $m_d(\mu)$ lies in $Z(\gS_d)$. Since $\mu \in \eP_d(l)$, it is clear that $m_d(\mu) \in \C_{CX}$, so we just need to check that it commutes with each $s_i$ for $i \in I_d$. Since $s_i$ commutes with $y_j(k)$ unless $j = i$ or $j=i+1$, we need only check that $s_i$ commutes with the elements $y_i (k) y_{i+1}(k)$ and $y_i (k) y_{i+1}(m) + y_i (m) y_{i+1} (k)$ for even values of $k$ and $m$. These calculations are similar: we consider the first. Write $k = (a-l)l + r$ as usual. Then
\begin{equation}\label{mbasis1} y_i(k) y_{i+1}(k)  =  \sum_{\begin{array}{c} i_1 \dots i_{a-1} < i \\ j_1,\dots,j_{a-1}< i+1 \end{array}} (i_1 \dots i_{a-1} \ i)^{(r)} (j_1 \dots \ j_{a-1} \ i+1 )^{(r)} \end{equation}
whereas
\begin{eqnarray*} s_i y_i(k) y_{i+1} s_i & = & \sum (i_1 \dots i_{a-1} \ i+1)^{(r)} s_i (j_1 \dots \ j_{a-1} \ i+1 )^{(r)} s_i \\
& = & \sum_{j_1,\dots,j_{a-1} \neq i} (i_1 \dots i_{a-1} \ i+1)^{(r)} (j_1 \dots \ j_{a-1} \ i )^{(r)} \\
& & + \sum_{j_1,\dots,j_{a-1} \neq i} (i_1 \dots i_{a-1} \ i+1)^{(r)} (j_1 \dots i+1 \dots \ j_{a-1} \ i )^{(r)}. \end{eqnarray*}
Showing that this equals (\ref{mbasis1}) is an exercise in applying Proposition~\ref{cxmult1}, Theorem~\ref{cxshift}, and (if $r=0$) Example~\ref{simpleoverlap2} according to whether 
$$\{i_1, \dots, i_{a-1} \} \cap \{j_1,\dots,j_{a-1}\}$$
equals $0$, $1$, or $2$. 

To complete the proof, we claim that $m_d(\mu) = z_d(\phi^{-1}(\mu)) + (*)$, where $(*)$ denotes a linear combination of $z_d(\nu)$s with strictly greater redundancy than $\phi^{-1}(\mu)$. By definition, $m_d(\mu)$ is a sum of products of CX-cycles of the form $A_1^{(r_1)} \dots A_a^{(r_a)}$. If the index sets are pairwise distinct, this contributes something of cycle type $\phi^{-1}(\mu)$, and certainly every element of cycle type $\phi^{-1}(\mu)$ arises in such a way, which provides the $z_d(\phi^{-1}(\mu)$. If the index sets overlap, any nonzero contribution that arises will be a similar product of CX-cycles such that $|A_1| + \dots +|A_a| < |\mu|$. Consequently, their cycle types will have strictly greater redundancy. \end{proof}

\section{Center of $\S_d^f$}\label{final}

\begin{Lemma}\label{JM}
Let $k = (a-1)l + r$ be even. Then $\x_i^k$ lies in $F_{(a-1)(l-1)+r}(\S^f_d)$. Moreover, we have
$$ \gr_{(a-1)(l-1)+r}(\x_i^k) = y_i(k) . $$
\end{Lemma}

\begin{proof}
This is the (slightly cleaner) analogue of Lemma 3.1 in \cite{B}. We use the same strategy, looking first at the case $k=(a-1)l$ and proceeding by induction on $a$. The base case $a=2$ is Lemma~\ref{lcase}.  

For the induction step, suppose first that $a=3$, that is that $k=2l$. Write $\x_i^{2l} = \x_i^l \x_i^l$ and apply Lemma~\ref{lcase} again. Unfortunately, if $l$ is odd then none of the terms that arise are CX-cycles, so Proposition~\ref{cxmult1} does not apply. However, for distinct $i,j,k \in I_d$ it is easy to check that
$$ (j \ i)^{(0)} (j \ k)^{(0)} = (i \ j \ k )^{(0)} ,$$
which is enough to establish that $\gr_{2(l-1)} \x_i^{2l}$ is of the correct form. 

For the case $a>3$ (and the case $a=3$ if $l$ is even) one can write $\x_i^k = \x_i^{(a-2)l} \x_i^{2l}$ and use Proposition~\ref{cxmult1}. 

Finally, if $k=(a-1)+r$ for $r > 0$, write $\x_i^k = \x_i^r \x_i^{(a-1)l} $ and apply the previous argument. \end{proof}

\begin{Definition}{\rm
For $\mu \in \eP_d(l)$, define the element $p_d(\mu)$ of $\S_d$ to be the symmetric polynomial
$$ \sum_{(\nu_1,\dots,\nu_d) \sim \mu} \x_1^{\nu_1} \dots \x_d^{\nu_d} \in Z(\S_d) .$$
}\end{Definition}

\begin{Lemma}\label{ptom}
For $\mu \in \eP_d(l)$, the image of $p_d(\mu)$ in $\S^f_d$ lies in  $F_{|\mu| - |\mu / l| } (\S^f_d)$, and 
$$ \gr_{|\mu| - |\mu / l|}(p_d(\mu)) = m_d(\mu) .$$     
\end{Lemma}

\begin{proof}
The polynomial degree of $p_d(\mu)$ in $\S_d$ is clearly $|\mu|$, and so the result follows from Lemma~\ref{JM}.
\end{proof}

\begin{Theorem} 
Let $l$ be odd. Then 
$$ \gr( Z(\S^f_d) ) = Z(\gS_d) ,$$
and the set $ B = \{ \ p_d(\mu) \ | \ \mu \in \eP_d(l) \ \}$ is a basis for $Z(\S^f_d)$. 
\end{Theorem}

\begin{proof}
Our filtration of $\S^f_d$ induces a filtration of $Z(\S^f_d)$ in which, for each $\mu \in \eP_d(l)$, the element $p_d(\mu)$ lies in degree $|\mu| - |\mu/l|$. Theorem~\ref{mbasis} and Lemma~\ref{ptom} establish that the set 
$$ \{\  \gr_{|mu|-|\mu/l|} (p_d(\mu)) \ | \ \mu \in \eP_d(l) \ \} $$
is a basis for the associated graded object $\gr( Z(\S^f_d))$. This is sufficient to establish that $B$ is a basis for $Z(\S^f_d)$ as claimed. \end{proof}

\end{document}